\documentclass[12pt, leqno]{article}
\usepackage{amsmath}
\usepackage{amsfonts}
\usepackage{amssymb}
\usepackage{graphicx}
\usepackage[monochrome]{color}
\usepackage{color}
\newtheorem{thm}{Theorem}[section]
\newtheorem{cor}[thm]{Corollary}
\newtheorem{lem}[thm]{Lemma}
\newtheorem{prop}[thm]{Proposition}

\numberwithin{equation}{section}

\newcommand{\Es}{{\mathcal E}}

\newcommand{\RR}{\mathbb{R}}
\newcommand{\ren}{\mathbb{R}^n}

\newcommand{\wu}{\widetilde{u}}
\renewcommand{\wp}{\widetilde{p}}
\newcommand{\ve}{\varepsilon}

\newcommand{\dint}{\displaystyle\int}

\def\qed{\,\unskip\kern 6pt \penalty 500
\raise -2pt\hbox{\vrule \vbox to8pt{\hrule width 6pt
\vfill\hrule}\vrule}\par}
\definecolor{darkblue}{rgb}{0.05, .05, .65}
\definecolor{darkgreen}{rgb}{0.05, .70, .05}
\definecolor{darkred}{rgb}{0.8,0,0}
\def\qed{\unskip\kern 6pt \penalty 500
\raise -2pt\hbox{\vrule \vbox to8pt{\hrule width 6pt
\vfill\hrule}\vrule}\par}
\begin{document}
\title{\textbf{Asymptotic behaviour of a porous \\ medium equation
with \\ fractional diffusion}}
\author{\Large Luis Caffarelli\footnote{caffarel@math.utexas.edu}
~and~ Juan Luis Vazquez\footnote{juanluis.vazquez@uam.es}\color{blue}\footnote{The main results of this paper have been announced at the International Conference ``Free Boundary Problems, Theory And Applications'' (FBP08), held in Stockholm in June 2008.}}

\date{ }
\maketitle

\begin{abstract}
We consider a porous medium equation with nonlocal diffusion effects
given by an inverse fractional Laplacian operator. In a previous paper we have found mass-preserving, nonnegative
weak solutions of the equation satisfying energy estimates. The equation is posed in the whole space $\ren$. Here we establish the large-time behaviour.
We first find selfsimilar nonnegative solutions by solving an elliptic obstacle problem with fractional
Laplacian for the pair pressure-density, which we call obstacle Barenblatt solutions. The theory for elliptic
fractional problems with obstacles has been recently established. We then use entropy methods to show
that the asymptotic behavior of general finite-mass solutions is described after renormalization by these special solutions,
which represent a surprising variation of the Barenblatt profiles of the standard porous medium model.

\end{abstract}

\newpage


\section{Introduction}
\label{sec.intro}

In the paper \cite{CaVa09} we have introduced a model of nonlinear diffusion with nonlocal effects  given by the system
\begin{equation}\label{eq1}
u_t=\nabla\cdot(u\nabla p), \quad p={\cal K}(u).
\end{equation}
Here, $u$ is a function of the variables $(x,t)$ to be thought of
as a density or concentration, and therefore nonnegative, while
$p$ is the pressure, which is related to $u$ via a linear operator
${\cal K}$, which we assume to be the inverse of a fractional Laplacian, ${\cal
K}=(-\Delta)^{-s}$, $0<s<1$.  The problem is posed for $x\in \RR^n$, $n\ge 1$, and $t>0$, and we
give initial conditions
\begin{equation}\label{eq.ic}
u(x,0)=u_0(x), \quad x\in \RR^n,
\end{equation}
where $u_0$ is a nonnegative and integrable function in $\RR^n$ with compact support
or fast decay at infinity. Motivation and related ideas for such a problem are given in \cite{CaVa09},
where the existence of a weak solution is established for a suitable class of initial data that includes all
functions $u_0$ with the above properties. Besides, a  number of the most important properties is proved, like the energy estimates, the bounds in the $L^p$ spaces and the property of finite propagation that says that compactly supported data produce solutions whose support is compact for all positive times. However, the  uniqueness of weak solutions is a pending open problem but for the case of one space dimension, \cite{BKM}. Comparison theorems, a crucial tool in parabolic
equations are only available under special circumstances (i.e., for so-called true super- or sub-solutions), as shown in \cite{CaVa09}. We will also need the results of \cite{CSV}, where we show that solutions with integrable data are actually bounded (so-called $L^1$-$L^\infty$ effect) and furthermore that bounded solutions are $C^\alpha $ continuous with H\"older exponent depending on the equation and local H\"older constant depending on the $L^\infty$ norm of the solution.

\medskip

In this paper we study  the asymptotic behaviour of such weak solutions. First, we introduce the rescaled flow and then we prove that the stable configurations of this modified equation are just the solutions of an {\sl obstacle problem} with fractional Laplacian, as studied by Athanasopoulos, Caffarelli, Salsa and Silvestre, \cite{ACS} and \cite{CSS}. This connection, described in Section \ref{sec.rf}, is quite surprising. In the original variables such equilibrium solution translates into  a selfsimilar solution that we propose to name ``fractional Barenblatt solution'' after the analogy with the porous medium case, though the analytical properties are quite different, see Theorem \ref{exist.obs}. In Sections \ref{sec.estimates2}, \ref{sec.ab} we prove  convergence of the rescaled flow to equilibrium by using a certain ``Boltzmann entropy'' technique; in the original variables this means that the asymptotic behaviour of finite mass, compactly supported solutions of our problem is given by the family of fractional Barenblatt solutions constructed from the  stationary fractional obstacle problem. Uniqueness of the asymptotic profiles is proved in Section \ref{sec.uniqueness}
and the proof of the main result, Theorem \ref{thm.asymp}, is completed in Section \ref{thm.asymp.end}. A final Section \ref{sec.sg} introduces the open question of spectral gap.

\medskip

\noindent {\bf Notation.}  We will use the notation $L_s=(-\Delta)^{s}$ with $0<s<1$ for the
fractional powers of the Laplace operator defined by Fourier transform  on the Schwartz class of smooth functions in $\RR^n$ and extended in a natural way to functions in the Sobolev space $H^{2s}(\RR^n)$.  Technical reasons of paper \cite{CaVa09} imply that in one space dimension the restriction $s<1/2$ will be observed. The inverse operator is denoted by ${\cal K}_s=(-\Delta)^{-s}$ and can be realized by convolution
$$
{\cal K}_s=K_s\star u, \qquad K_s(x)=c(n,s)|x|^{2s-n}.
$$
see  \cite{Stein}. ${\cal K}_s$ is a positive self-adjoint operator. We will write ${\cal H}_s={\cal K}_s^{1/2}$ which has kernel $K_{s/2}$. The subscript $s$ will be omitted when $s$ is fixed and known. For functions that depend on $x$ and $t$, convolution is applied for every fixed $t$ with respect to the space variables. We then use the abbreviated notation $u(t)=u(\cdot,t)$.
\normalcolor

\medskip

\color{blue} \noindent {\sc Related works.} Some papers have  appeared since 2008 dealing with similar equations. As closely related, we mention the recent research note by Biler et al. \cite{BIK} where a more general version of our model is proposed, and existence of solutions, decay estimates and selfsimilar solutions are announced. Full details are still to be provided. At least in the latter topic, their method is different and interesting. The asymptotic behaviour is not tackled. \normalcolor

\section{Existence result and basic estimates}
\label{sec.estimates1}

\noindent {\bf Definition.} We say that $u$ is a weak solution of equation (\ref{eq1})
in $Q_T=\RR^n\times (0,T)$ with initial data $u_0\in L^1(\RR^n)$ if $u\in L^1(Q_T)$, ${\cal K}(u)\in W^{1,1}_{loc}(Q_T)$, and $u\,\nabla{\cal K}(u)\in L^1(Q_T)$,
and the identity
\begin{equation}
\iint u\,(\phi_t-\nabla {\cal K}(u)\cdot\nabla\phi)\,dxdt+ \int
u_0(x)\,\phi(x,0)\,dx=0
\end{equation}
holds for all continuous test
functions $\phi$ in $Q_T$ such that $\nabla_x\phi$ is continuous, and $\phi$ has compact support in the space variable and vanishes near $t=T$. \color{blue} The following result is proved in \cite{CaVa09}. When $T=\infty$ we write $Q$ instead of $Q_T$.

\begin{prop}\label{thm.ex}  Let $u_0\in L^1(\ren)\cap L^\infty(\ren)$, $u_0\ge 0$, and such that
\begin{equation}
u_0(x)\le A\,e^{-a|x|} \qquad \mbox{for some $A,a>0$}\,.
\end{equation}
Then there exists a weak solution $u$ of Equation \eqref{eq1} with initial data $u_0$. Besides,  $u\in L^\infty(0,\infty: L^1(\ren))$, $u\in L^\infty(Q)$, $\nabla {\cal H}(u)\in L^2(Q)$. For all $t>0$ we have
\begin{equation}
\int_{\ren} u(x,t)\,dx=\int_{\ren} u_0(x)\,dx\,,
\end{equation}
and $\|u(t)\|_\infty \le \|u_0\|_\infty$. The solution decays exponentially as $|x|\to\infty$. The first energy inequality holds in the form
\begin{equation}
\dint_0^{t_1}\dint_{\ren} |\nabla {\cal  H} u|^2\,dxdt +  \dint_{\ren}  u(t_1)\log(u(t_1))\,dx
 \le \dint_{\ren} u_0\log(u_0)\,dx\,,
\end{equation}
while the second says that for all $0<t_1<t_2<\infty$
\begin{equation}
\int_{t_1}^{t_2}\int_{\ren} u\,|\nabla {\cal K}u|^2\,dxdt+ \frac 12\int_{\ren} |{\cal H}u(t_2)|^2\,dx\le
\frac 12\int_{\ren} |{\cal H}(u(t_1)|^2\,dx\,.
\end{equation}
\end{prop}

\noindent {\bf Properties of the constructed solutions.} Here are some of the most useful

 \noindent - Conservation of mass
\begin{equation}
\frac{d}{dt}\int u(x,t)\,dx=0.
\end{equation}

\medskip

\noindent - $L^p$ estimates. We also prove that the
$L^p$ norm of the solutions, $1< p\le \infty$, does not increase in time.

\medskip

\noindent - Conservation of sign: $u_0\ge 0$ implies that $u(t)\ge 0$
for all times.

\medskip

\noindent - Persistence of strict positivity: At every point $x_0$ where  $u(x_0,t_0)>0$ for some time $t_0$
we have $u(x_0,t)>0$ for all later times $t>t_0$.

\medskip

\noindent - {Finite propagation. Solutions with compact support:}
One of the most important features of the porous medium equation
and other related degenerate parabolic equations is the property
of finite propagation, whereby  compactly supported initial data
$u_0(x)$ give rise to solutions $u(x,t)$ that have the same
property for all positive times, i.e., the support of $u(\cdot,t)$
is contained in a ball $B_{R(t)}(0)$ for all $t>0$.

\begin{prop}\label{prop.sc1} {\rm \cite{CaVa09}} Assume that $u$ is a bounded solution, $0\le u\le L$, of equation \eqref{eq1} with \  ${\cal K}=(-\Delta)^{-s}$ with $0<s<1$ $(0<s<1/2$ if $n=1)$, as constructed in Theorem {\rm \ref{thm.ex}}.  Assume that $u_0$   has compact support. Then $u(\cdot,t)$ is compactly supported for all $t>0$. More precisely, if  \ $0<s<1/2$ and $u_0$ is below the ''parabola-like'' function
\begin{equation}
U_0(x)=a(|x|-b)^2,
\end{equation}
for some $a,b>0$, with support in the ball $B_b(0)$, then there is $C(n,s,L,a,T)$  large enough, such that
\begin{equation}
u(x,t)\le a(Ct-(|x|-b))^2
\end{equation}
\color{blue} for $x\in \RR$ and $0<t<T$. \normalcolor
\end{prop}

 \medskip

\noindent  - A standard comparison result for parabolic equations does not seem to work. This is one of the main technical difficulties in the study of this equation. In fact, we find special situations where some comparison holds by using so-called true super- and subsolutions.

 \medskip

 \noindent  - Next, we state the $L^1$ to $L^\infty$ result that is  proved in the forthcoming paper \cite{CSV}.

\begin{prop}\label{th:L-inf} Let $u$ be a weak solution of Problem \eqref{eq1}--\eqref{eq.ic} with $u_0\in L^1(\mathbb{R}^N)\cap L^\infty(\mathbb{R}^N)$, and $u_0$ decreases exponentially as $x\to\infty$. Then  there exists a positive constant $C$ such that for every $t>0$
\begin{equation}
\sup_{x\in\mathbb{R}^N}|u(x,t)|\le C\,t^{-\alpha }\|u_0\|_{L^1(\mathbb{R}^N)}^{\gamma}
\label{eq:L-inf}\end{equation}
with $\alpha=n/(n+2-2s)$, $\gamma=(2-2s)/((n+2-2s)$. The constant $C$ depends only on $n$ and $s$.
\end{prop}

\normalcolor

\noindent - {\bf Regularity.} Bounded solutions are $C^\alpha$ smooth for some $\alpha>0$: The H\"older exponent depends on $n$ and $s$ and the local H\"older constant depends also on the $L^\infty$ norm of the solution. This is proved in \cite{CSV}.

 \medskip

\newpage

\section{Large time behavior. Rescaled flow}
\label{sec.rf}

\color{blue}

We begin here the study of the large time behavior. As a first step in the analysis we will introduce self-similar variables, typical in porous medium theory, leading to a rescaled evolution equation. In the analysis of the steady states of that flow we will discover the solutions of an elliptic obstacle problem with fractional diffusion. This surprising connection is the main novelty of the paper, and the basis for the analysis of stabilization to be done in subsequent sections. The precise obstacle result is carefully stated in Theorem \ref{exist.obs}.
\normalcolor

We take a weak solution $u\ge 0$ with integrable and bounded initial data, as constructed in \cite{CaVa09}. Inspired by the asymptotics of the standard porous medium equation, we define the {\sl rescaled flow} through the transformation
\begin{equation}
u(x,t)=(1+t)^{-\alpha} v(x\,(1+t)^{-\beta}, \tau)
\end{equation}
with new time $\tau=\log(1+t)$. We also put $y=x\,(1+t)^{-\beta}$ as rescaled space variable.
In order to cancel the factors including $t$ in a explicit way,
we  get the condition on the exponents
\begin{equation}
\alpha + (2-2s)\beta =1,
\end{equation}
where we use the homogeneity of ${\cal K}$ in the form
\begin{equation}
({\cal K} u)(x,t) = (1+t)^{-\alpha+2s\beta} ({\cal K}v)(y,\tau).
\end{equation}
Since we also want conservation of (finite) mass, we must put
$\alpha= n \beta$. We get the precise value of the exponents:
\begin{equation}\label{def.a.b}
\beta=1/(n+ 2-2s), \quad
\alpha=n/(n+ 2-2s).
\end{equation}
We recall that $0<s<1$ so that $\beta\in (1/(n+2), 1/n)$. The first value is the standard porous medium case.
In this way,  we arrive at the {\sl nonlinear, nonlocal
Fokker-Planck equation}
\begin{equation}\label{ren.eq}
v_\tau=\nabla_y\cdot(v\,(\nabla_y {\cal K}(v)+\beta y ))
\end{equation}
with $\beta$ as given above.  Note that this formula implies a
transformation for the pressure of the form
\begin{equation}
p(u)(x,t)=(1+t)^{-\sigma} p(v)(x\,(1+t)^{-\beta}, \tau), \quad \mbox{with }
\sigma=\alpha-2s\beta=1-2\beta=\frac{n-2s}{n+2-2s}<1.
\end{equation}
We are not making much use of this last formula.

\medskip

\noindent $\bullet$ {\bf Equilibrium states for the rescaled flow}

We want to find stationary solutions of the
rescaled equation \eqref{ren.eq}, i.e., solutions $V(y)$ of
the system
\begin{equation}\label{baren.eq}
\nabla_y\cdot(V\,\nabla_y (P+ a |y|^2))=0, \quad P={\cal K}(V).
\end{equation}
where $a=\beta/2$, and $\beta$ defined just above. Since we are
looking for asymptotic profiles of the standard solutions of
(\ref{eq1}) we also want $V\ge 0$ and integrable. The simplest
possibility is integrating once and getting the radial version
\begin{equation}\label{baren.eq2}
V\,\nabla_y (P+  a |y|^2))=0, \quad P={\cal K}(V), \quad V\ge 0.
\end{equation}
The first equation gives an alternative choice that reminds us of the
complementary formulation of the obstacle problems, \cite{FdmObst}, \cite{CaffObst}.

\medskip

\noindent $\bullet$ {\bf  Obstacle problem.
Barenblatt solutions of new type}
 \label{ssec.eq}

Indeed, if we solve the obstacle problem with
fractional Laplacian we will obtain a unique solution $P(y)$ of
the problem:
\begin{equation}\label{obst}
\begin{array}{l}
P\ge \Phi, \quad V=(-\Delta)^{s} P\ge 0; \\
\mbox{either } \ P=\Phi \ \mbox{or } \ V=0.
\end{array}
\end{equation}
with $0<s<1$.  In our present application we have to choose as obstacle $\Phi=C- a |y|^2$, where $C$ is any positive constant and $a=\beta/2$. For uniqueness we also need the condition $P\to 0$ as $|y|\to\infty$.
 The theory is developed in the papers \cite{ACS, CSS}, the solution of this obstacle problem is unique and belongs to the space $H^{-s}$ with pressure in $H^{s}$. The  solutions have the following regularity:  $P\in C^{1,s}(\ren)$ and $V\in C^{1-s}(\ren)$. The figure represents an approximate plot of $P(y)$.

\centerline{
{\includegraphics[width=9.5 cm,height=5.5 cm]{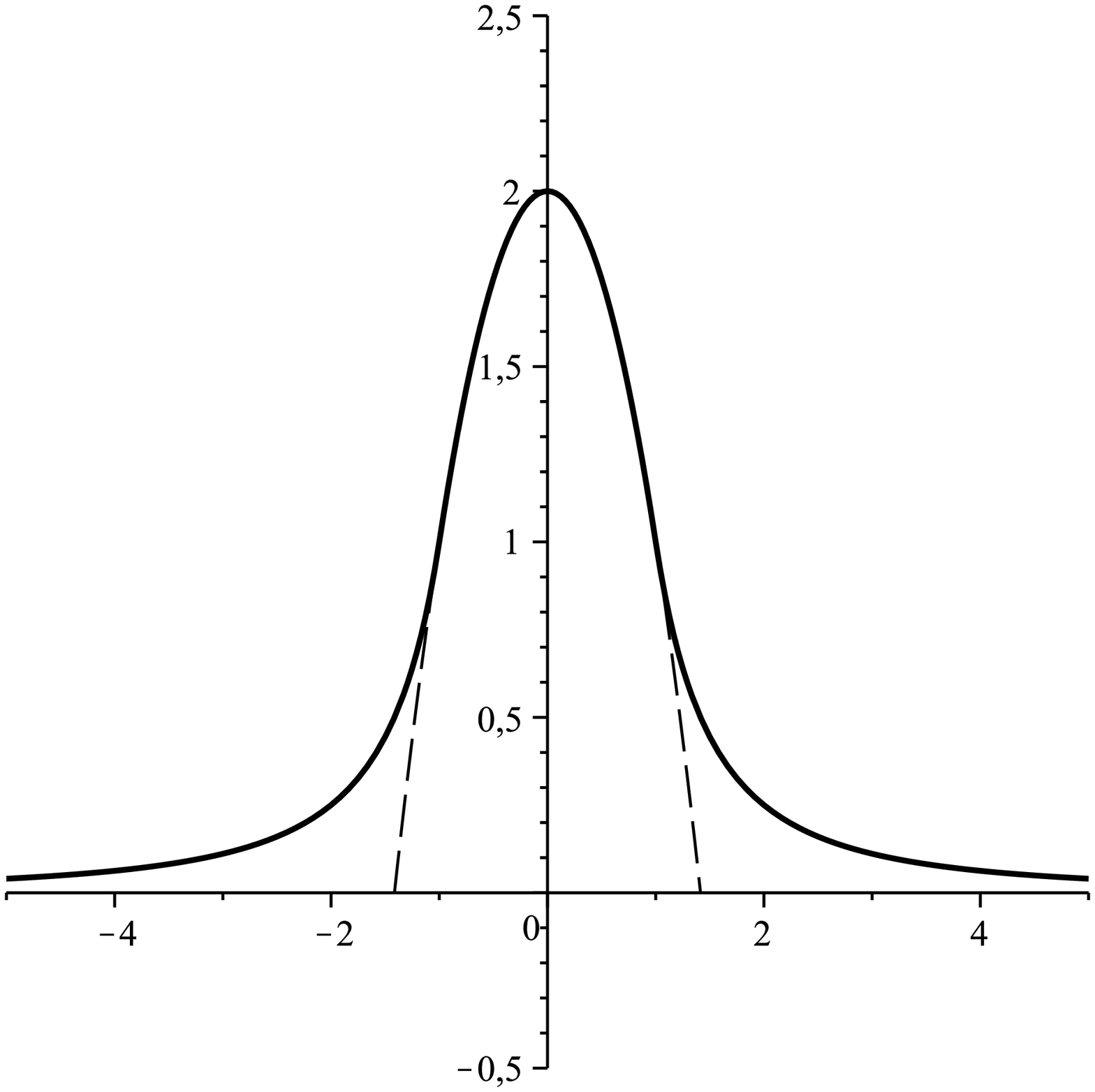}}
}

Note that for $C\le 0$ the solution is trivial, $P=0$, $V=0$,
hence we choose $C>0$. We also note that the pressure is defined
but for a constant, so that we may take without loss of generality
$C=0$ and take as pressure $\widehat P= P-C$ instead of $P$. But
then $P\to 0$ implies that $\widehat P\to -C$ as $|y|\to\infty$,
so we get a one parameter family of stationary profiles that we
denote $V_C(y)$.  In any case, if we consider the free boundary
location, $|x|=R(C,\sigma,n)$, which is the boundary of the
contact set of the obstacle problem, ${\cal C}=\{|y|\le R(C)\}$,
then $R$ tends to $0$ as $C\to 0$, and $R$ tends to $\infty$ as
$C\to\infty$. Let us summarize the results of this section.

\begin{thm}\label{exist.obs} For every $C>0$ there exists a unique solution $P=P_C(y)\in C^{1,s}$ of the obstacle problem \eqref{obst}
with $\Phi(y)=C-a \,|y|^2$ and $P\to0 $ as $|y|\to\infty$.  It is radially symmetric as a function of $y$ and we have the scaling law
\begin{equation}
P_C(y)=C\,P_1(y/C^{1/2}).
\end{equation}
Moreover, $V_C=(-\Delta)^s P_C$ satisfies the scaling relation $V_C(y)=C^{1-s}V_1(y/C^{1/2})$ and $V_1(y)$ has compact support. Finally, if we put $a=\beta/2>0$, then
\begin{equation}
U_C(x,t)= (1+t)^{-\alpha} V_C(x\,(1+t)^{-\beta})
\end{equation}
is a weak solution of Equation {\rm (\ref{eq1})}.
\end{thm}

Note that  $U_C(x,t)$ is a weak solution of Equation (\ref{eq1})
in the sense we have defined. It lives for $t>-1$ and it takes at the initial time $t=-1$
a Dirac delta $M\,\delta(x)$ as initial data, i.e., it is the
source-type or Barenblatt solution for this problem, and it has a
profile $V_C\ge 0$ that has compact support. In the sequel we call these solutions  the family of
{\sl  fractional Barenblatt solutions.}

Some further observations: an easy calculation gives the realationship
\begin{equation}\label{mass}
M=\int U_C(x,t)\,dx=\color{blue} \int V_C(y,\tau)\,dy\normalcolor = c(n,s)\,C^{(n+2-2s)/2}
\end{equation}
(compare with the PME case in \cite{JLVSmoothing}, page 23). Clearly, $V\ge 0$  is supported in the contact set of the
obstacle problem, ${\cal C}=\{|y|\le R(C)\}$. The radius $R$ is smaller than the intersection of the parabola $\Phi$ with the axis
$R_1=(C/a)^{1/2}$. On the other hand, the pressure $P_C(|y|)$ is
always positive and decays to zero as $|y|\to\infty$ according to
fractional potential theory, \cite{Land, Stein}.


\section{Entropy estimates for the rescaled problem}
\label{sec.estimates2}

\normalcolor

 Our next problem is now to prove that these profiles are attractors for the rescaled flow that we have introduced. In the next sections we prove  the main result of this paper.

 \begin{thm}\label{thm.asymp} Let $u(x,t)\ge 0$ be a weak solution of Problem \eqref{eq1}--\eqref{eq.ic} with bounded and integrable initial data such that $u_0 \ge 0$ has finite entropy in the sense defined in formula \eqref{entropyE}. Let $v(y,\tau)$ be the corresponding rescaled solution. As $\tau\to\infty$ we have
 \begin{equation}
 v(\cdot,\tau)\to V_C(y) \quad \mbox{in \  $ L^1(\ren)$ and also in $L^\infty(\ren)$}.
 \end{equation}
 The constant $C$ is determined by the rule of mass equality: $\int_{\ren} v(y,\tau)\,dy=\int _{\ren} V_C(y)\,dy$.
 In terms of $u$ this translates into
 \begin{equation}
 u(x,t)-U_C(x,t)\to 0 \quad \mbox{in } L^1(\ren), \quad t^{\alpha}|u(x,t)-U_C(x,t)|\to 0 \quad \mbox{uniformly in $x$},
 \end{equation}
 both limits taken as $t\to \infty$.
 \end{thm}

As a first step in the proof, we review the estimates of Section \ref{sec.estimates1} in order to adapt them to the rescaled problem. There is no problem is
reproving mass conservation or positivity for the rescaled flow, we leave it to the reader checking both facts.
By shifting a bit the origin of time, there is no loss of generality in assuming that the entropies mentioned below are finite even at $t=0$ even if they were not assumed to be so. The values of $\alpha$ and $\beta$ are those fixed in \eqref{def.a.b}.

\noindent $\bullet$ The first energy estimate becomes  (recall the notation ${\cal H}={\cal K}^{1/2}$)
\begin{equation}
\begin{array}{c}
\displaystyle \frac{d}{d\tau }\int v(y,\tau )\log v(y,\tau )\,dy=
-\int
|\nabla {\cal H}v|^2\,dy-\beta \int \nabla v\cdot y \\
\displaystyle = -\int |\nabla {\cal H}v|^2\,dy + \alpha \int v.
\end{array}
\end{equation}

\noindent $\bullet$ However, the second energy estimate has an
essential change. We need to define the entropy of the rescaled flow as
\begin{equation}\label{entropyE}
\Es(v(\tau)):=\frac{1}2\int_{\RR^n} (v\,{\cal K}(v) + \beta |y|^2v)\,dy
\end{equation}
The entropy contains two terms. The first is
$$
\Es_{1}(v(\tau )):=\int_{\RR^n} v\,{\cal K}(v)\,dy=\int_{\RR^n} |{\cal H}v|^2\,dy,
$$
hence  positive. The second is the moment $\Es_{2}(v(\tau
))=M_2(v(\tau )):=\int |y|^2v\,dy$, also positive.
 By differentiation we get
\begin{equation}
\frac{d}{d\tau }\Es(v)=-{\mathcal I}(v), \quad \quad {\mathcal I}(v):=\int
\left|\nabla ({\cal K}v + \frac{\beta}2 |y|^2)\right|^2\,vdy.
\end{equation}
We call ${\mathcal I}(v)$ the entropy dissipation of the rescaled solution $v$.
This is formal, the rigorous formula says that whenever the initial entropy is finite, then
$\Es(v(\tau ))$ is uniformly bounded for all $\tau >0$, ${\mathcal I}(v)$
is integrable in $(0,\infty)$, and
\begin{equation}
\Es(v(\tau ))+\int_0^\tau\int_{\RR^n}\left|\nabla ({\cal K}v + \frac{\beta}2
|y|^2)\right|^2\,vdy\,dt\le \Es(v_0).
\end{equation}
The bound on $\Es(v(\tau ))$ implies a bound for the energy
${\cal E}_{1}(v(\tau ))$ and for the moment $M_2(v(\tau ))$ in
$L^\infty(0,\infty)$. We also conclude that
\begin{equation}
\int_0^\infty d\tau \int_{\RR^n} dy \,v\,| \nabla {\cal ({\cal K}}v + (\beta/2)
|y|^2) |^2 <\infty.
\end{equation}
 This estimate shows that the integral \ ${\mathcal I}(v(\tau)) $ converges to zero as $\tau\to 0$ in some kind of time average, which is a basic fact in the asymptotic analysis.

\medskip

\noindent {\bf Other consequences:} (1)  For the original variable: the
fact that $\int v(y,\tau) |y|^2dy$ is bounded uniformly in time implies in terms of the original variable $u$ that
$$
\int u(x,t)\,x^2dx \le C\,t^{2\beta}.
$$

\noindent (2) The fact that $\int |{\cal H}(v)|^2dy$ is bounded uniformly in time
implies that
$$
\int |{\cal H}(u)|^2dy \le C\,t^{2s\beta-\alpha}.
$$

\section{Asymptotic behavior. Zero entropy dissipation}
\label{sec.ab}
\normalcolor

We take a solution $u(x,t)$ under the same assumptions on the initial data, we rescale it as in the previous sections to obtain $v(y,\tau)$, and we put $\int v_0(y)\,dy=M>0$.  Now the idea is to let $\tau \to \infty$ in the rescaled flow equation \eqref{ren.eq}. Since the entropy \eqref{entropyE} is nonnegative and goes down, there is a limit
$$
E_*=\lim_{t\to\infty} {\cal E}(\tau)\ge 0.
$$
Moreover,  $v$ is bounded in $L_y^1$ uniformly in $\tau$,  and  also $\int |\nabla {\cal H}(v)|^2\,dy$ is  uniformly bounded, in time, therefore we have that $v(\cdot,\tau)$ is a compact family and there is a subsequence $v(\cdot,\tau_j)$  with $\tau_j\to\infty$ that converges in $L_y^1$ and almost everywhere to a limit $v_*\ge 0$.

In this section we prove the following intermediate result:

\begin{lem} the mass of $v_*$ is the same as the mass of $v(\cdot,\tau)$, $\int v_*(y)\,dy=M$ (hence, in particular, $v_*$ is not trivial). The functions $v_*(y,\tau)$ and $w(y,\tau) ={\cal K}v_*+ \frac{\beta}2 |y|^2$ are continuous and $w(\cdot,\tau)$ is constant on every connected component of the set where $v_*$ is not zero,
\end{lem}

\noindent {\sl Proof.}  We want to use compactness of the family of time-translates
$$
v_j(y,\tau)=v(y,\tau+\tau_j), \qquad \tau_j\to\infty.
 $$
 We observe that  the $L^\infty$ estimate of Theorem \eqref{th:L-inf} with sharp exponents implies that $v(y,\tau)$ is uniformly bounded in time (for all  times $\tau\ge 1$), hence the family $v_j$ is uniformly bounded. Moreover, bounded families of solutions are uniformly equicontinuous (H\"older continuous), and this applies both to $u$ and $v$.
 This was proved as part of the local regularity of paper \cite{CSV}. Therefore, the convergence $v_j\to v_*$ takes place locally uniformly in $Q=\RR^n\times (0,\infty)$.

Next, we prove the conservation of mass. This a consequence of the uniform estimate of $v(\cdot,\tau)$ in $L^\infty(\ren)$  for all large $t$, the local compactness of the solutions, and finally  the uniform boundedness of the moment  $\int v(y,\tau)\,|y|^2 dy$. Another immediate consequence is that the {\sl lim inf} of the entropy  ${\cal E}(v(\tau_j))$ is equal or larger that $M_2(v_*)$.

We also have ${\cal H}(v)\in L^2_y$ uniformly in $\tau$ as well as
$\nabla {\cal H}(v)$ in $L^2_{y,\tau}$. Since ${\cal K}(v)={\cal H}({\cal H}(v))$ and $\nabla {\cal K}(v)=\nabla {\cal H}({\cal H}(v))={\cal H}(\nabla ({\cal H} v))$, we derive from the bound for $\nabla {\cal H}v$ in $L^2_{y,\tau}$  estimates for ${\cal K}(v)$. We recall that $\nabla {\cal H}(v)$ is a ''derivative of order $1-s$ of $v$'', and since $v$ is bounded, $v\in L^\infty_{y,\tau}$, we conclude that $v\in L^2_t H^{1-s}_{x,loc}$. By potential theory, it is then clear that ${\cal K}(v)\in L^2_tH_x^{1+s}$.

All of this can be used in passing to the limit in the term $\iint v (\nabla {\cal K}v)\nabla \phi\,dxdt$ as follows: we have the convergence of $v_{j}$ in $C([0,T]:L^2(B_R))$ together with the weak convergence of $p_j ={\cal K}v_{j}$ and $\nabla p_{j}$. In this way, we find that $v_*$ is a weak solution of the rescaled equation \eqref{ren.eq}.

Finally, we want to prove that this solution is stationary. For that we pass to the limit in the entropy dissipation bound that says that  for every $h>0$ fixed
$$
\lim_{j\to\infty} \int_{0}^{h}\int\left|\nabla ({\cal K}v_j(y,\tau) + \frac{\beta}2
|y|^2)\right|^2\,v_j(y,\tau) \, dy\,d\tau =  0.
$$
By the strong and weak convergence results that we have, it follows that
\begin{equation}
{\mathcal I}(v_*)= \int_{0}^{h}\int\big|\nabla ({\cal K}v_*(y,\tau) + \frac{\beta}2
|y|^2)\big|^2\,v_*(y,\tau) , dy\,d\tau =  0.
\end{equation}
This implies that if $w(y,\tau) ={\cal K}v_*+ \frac{\beta}2 |y|^2$ \ then \ $v_*|\nabla w|^2=0$ \ a.\,e. in $Q$. It follows that for almost every time the continuous function $w$ must be a constant in space in every connected component of the set where the continuous function $v_*$ is not zero, $w(\cdot,\tau)=K(\tau)$. By continuity of $w$, this happens indeed for all times, but the space constant might still vary continuously in time.

The determination of the space constant, its independence of time, and finally the precise form of $v_*$ will allow us the unique identification of the limit as a solution of the Barenblatt obstacle problem with the same mass $\int v_*(y)\,dy=M$. This is done in the next section. and the proof of the asymptotic theorem is completed in Section \ref{thm.asymp.end}.


\section{Uniqueness for the asymptotic limit}
\label{sec.uniqueness}

We will prove here that the vanishing of the entropy dissipation ${\mathcal I} $ characterizes the stationary solutions of the rescaled flow. We use the notations $u$ and $x$ instead of $v$ and $y$ in this ``elliptic section''. We also put $\beta=1$ for simplicity. Thus, we assume that $\wu(x)\ge 0$ \ is continuous in $\RR^n$, that $\wp=(-\Delta)^{-s}\wu$ is also continuous, and that we have the basic identity
$$
\int |\nabla (\wp-Q)|^2\wu\,dx=0,
$$
where  $Q(x)=-|x|^2/2$. Then

\begin{thm}  There is a $c_0>0$ such that $\wp$ equals $P(c_0)$, the solution of the fractional  obstacle problem with obstacle $\Phi=Q+c_0$. Moreover, $P(c_0)\to 0$ as $x\to\infty$.
\end{thm}

We proceed via a series of lemmas. We write $P_0=P(c_0)$.

\begin{lem}
$\wp$ attains its maximum only at $x=0$
\end{lem}

\noindent  {\sl Proof. } By the maximum principle, $\wp $ cannot  attain its maximum at a point where $\wu=0$, since then $(-\Delta)^s\wp=0$ there, and this is impossible at a maximum of $\wp$ according to the integral formula that describes the operator $(-\Delta)^s$ (unless $\wp$ is constant, which is not the case). Hence, $\wp$ must attain its maximum at a point where $\wu>0$. Besides, at the maximum $x_0$ we have $\nabla \wp(x_0)=0$, hence $\nabla Q(x_0)=0$, and this means that $x_0=0$.

\begin{cor}
(i) We have $\wu(0)>0$. (ii) $\wp=Q+c_0$ in a neighborhood $\Omega$ of $x=0$.
\end{cor}

\begin{lem} Let $P_0$ is the solution of the fractional obstacle problem for data $\Phi=Q+c_0$, then
 $P_0-Q$ is strictly convex.
\end{lem}

\noindent  {\sl Proof. } Let $e$ be a unit vector and let $h>0$. We have
\begin{eqnarray*}
&\frac{1}{2}\left[P_0(x+he)+P_0(x-he)\right]\ge \frac{1}{2}\big[2c_0-\frac12\left(|x+he|^2+|x-he|^2\right)\big] =\\
&c_0-\frac{1}{2}|x|^2-\frac{1}{2}h^2= Q(x)+c_0-\frac{1}{2}h^2.
\end{eqnarray*}
Therefore, since $P_0$ is a supersolution of $(-\Delta)^s p(x)\ge 0$, also the function
$$
\gamma(x,h)= \frac{1}{2}\left[P_0(x+he)+P_0(x-he)\right]+ \frac{1}{2}h^2
$$
is a supersolution,  and $\gamma(x,h)\ge Q(x)+c_0$. Now, the solution of the obstacle problem has the property of being the minimal supersolution lying above the obstacle. Hence, $\gamma (x,h)\ge P_0(x)$. Dividing by $h^2$  we get $D_{ee}P_0\ge -1$ for any direction $e$. Now, $D_{ee} P_0$ is a supersolution outside the coincidence set, so
$D_{ee}P_0  > -1$ there. \qed

\begin{lem}  On each component of the set $\{\wu>0\}$ we have
$$
\wp=Q+c
$$
for some constant (depending possibly on the component). In fact, it holds on the closure of each component. We also have $\wp\to 0$ as $|x|\to\infty$.
\end{lem}

\color{blue} The proof is easy using the previous arguments. \normalcolor

\begin{lem}
If $P_0$ is the corresponding solution of the obstacle problem for data \ $\Phi=Q+c_0$,
then $P_0\ge \wp$ and $P_0=\wp$ in the coincidence set of $P_0$.
\end{lem}

\noindent  {\sl Proof. }Let $P_t=P_0+t$ so that $P_t\to t>0$ as $|x|\to\infty$. Then, $P_t-Q$ is strictly convex outside the coincidence set (see proof above), thus the gradient never vanishes. The obstacle is now $\Phi_t=\Phi+t$ so there is
no coincidence between $P_t$ and $\Phi=Q(x)+c_0$. Let $\Lambda_0$  be the coincidence set of $P_0$, the closure of a ball $B_\rho(0)$, and let  $\Theta_0$ be the connected component of $\{\wp=Q+c_0\}$ near the origin.

From the previous estimates, for $t$ large  $P_t$ does not touch $\wp$. Let us now decrease $t$ and suppose that there is a first $t_0>0$ such that $P_{t_0}$ touches $\wp$.

The  possibility of first contact $P_t=\wp$ \ at a point $x_0$ such that $\wu(x_0)=0$ contradicts the maximum principle, unless $P_{t_0}$ and $\wp$ coincide and then $t_0=0$.

The first contact cannot happen either at a point $x_0$ where $(-\Delta)^s P_0=0$ and $\wu>0$ since then $\wp=Q+c_1$ ($c_1>c_0$) in a neighborhood of $x_0$, and we have already proved that $\nabla P_{t_0}$ cannot equal $\nabla Q$.

Finally, we have to consider the case $(-\Delta)^s P_0>0$ and $\wu>0$, where $P_t=Q+c_0+t$ while $\wp=Q+c$. This forces again $t_0=0$.  This concludes that proof that $P_0\ge \wp$.

In order to check that $\Lambda_0\subset \Theta_0$ we consider for contradiction the existence of a point $x_0\in B_\rho(0)$ such that $x_0\in \partial \Theta_0$. In that case $P_0-\wp$ has a zero minimum there with $(-\Delta)^s P_0(x_0)>0$, $(-\Delta)^s \wp)(x_0)=0$, hence $(-\Delta)^s(P_0-\wp)(x_0)>0$ which contradicts the maximum principle.
\qed

\begin{lem}
We have $\wp\ge Q+c_0$.
\end{lem}

\noindent  {\sl Proof. } Suppose that $\wp-Q-c_0$ has a negative minimum at some point $x_0$. We can translate $\wp$ \ up and bring it down until it touches $P_0$ by above. This must happen on the set where $u_0=0$, since $\Lambda_0\subset \Theta_0$. But this contradicts the maximum principle since $\wu\ge 0$ at $x_0$ and $P_0$ and $\wp$ are not equal.

Finally, we prove that $\wp$ is a supersolution and $\wp\ge P_0$ (the least supersolution). \qed


\section{End of proof of Theorem \ref{thm.asymp}}
\label{thm.asymp.end}

We now go back to the end of Section \ref{sec.ab} and revert to the notations $y$ and $v$. Once we have uniquely identified the pressure $p_*$ at a.e. time as the solution $P_0$ of the obstacle problem for a certain constant $c_0$ that might depend on time, the scaling laws (cf. \eqref{mass}) allow to uniquely identify $c_0$ as the constant
such that $\int (-\Delta^sP_0)(y)\,dy=M$. Hence, $c_0$ is fixed in time and we conclude by continuity that $p_*(y,\tau)=P_{c_0}(y)$.

 The above argument also identifies in a unique way the constant  asymptotic limit $v_*(y,\tau)=V_C(y)$. It follows that the convergence  $v_j\to V_c$ as $\tau_j\to\infty$ takes place independently of the subsequence $\tau_j\to \infty$, and  we conclude that $v_j(\cdot,\tau)\to V_C$ as $\tau\to\infty$. Since the families $v_j$ are compact in $L^\infty_{loc}(\ren)$ the convergence is uniform in that space. On the other hand, from the conservation of mass and the uniform bound for $\int |y|^2 v(y,\tau)\,dy$ we conclude that the convergence takes place also in $L^1(\ren)$.

We still need to check the uniform convergence of $v(\cdot,\tau)$ to $V_C$ in the whole space. This is a rather simple calculus lemma using the following facts:

(i) for $\tau\ge 1$ the funtions $v(\cdot,\tau)$ are uniformly bounded,

(ii) for every $\ve>0$ there exists $R=R(\ve)$ such that  $\int_{|y|\ge R(\ve)} v(y,\tau)\,dy<\ve$ uniformly in $\tau$.

(iii) uniformly bounded functions are $C^{\alpha}$ continuous with a uniform bound.

This ends the proof for $v$. The translation of the results in terms of the variable $u$ is immediate.


\section{The question of a spectral gap}
\label{sec.sg}

When using the entropy--entropy dissipation approach, we are interested in the following inequality
\begin{equation}
\Es(v)\le C \,{\mathcal I}(v)
\end{equation}
for a class $\cal C$ of functions that includes our solutions and for a constant $C>0$ which is called the spectral gap. In our case it means that for all $v\in {\cal C}$ we must have
\begin{equation}
\int_{\RR^n} |{\cal H}v|^2\,dy+\beta \int_{\RR^n}  |y|^2v\,dy\le  C \int
\left|\nabla ({\cal K}v + \frac{\beta}2 |y|^2)\right|^2\,vdy.
\end{equation}
Developing the last integral leads to
\begin{equation*}
\int |{\cal H}v|^2\,dy+\beta \int  |y|^2v\,dy\le  C \int
\left|\nabla {\cal K}v\right|^2\,vdy + C\beta^2 \int |y|^2vdy + 2C\beta\int (\nabla {\cal K}v.y)vdy.
\end{equation*}
We can transform the last integral into
\begin{equation*}
\int (\nabla {\cal K}v.y)\,vdy=-n\int {\cal K}v\,vdy-\int  ({\cal K}v)(\nabla v. y)\,dy.
\end{equation*}
For many solutions the last term will be positive, hence we need
\begin{equation*}
(1+2Cn\beta) \int |{\cal H}v|^2\,dy \le  C \int
\left|\nabla {\cal K}v\right|^2\,vdy + (C\beta^2-\beta) \int |y|^2vdy.
\end{equation*}
The question is: can we prove this strange interpolation inequality? Since the order of differentiation in the first
term is $-s$, for the second $1-2s$, while the third is zero. It is clear that the case $s=1/2$ is different and maybe simpler.

\medskip


\

\noindent \textsc{Acknowledgment.} This work was initiated while the second author was a visitor
at ICES, Univ. of Texas, as an Oden Fellow. The second author supported by Spanish
Project MTM2008-06326-C02-01  and by ESF Programme ``Global and
geometric aspects of nonlinear partial differential equations".

\vskip 1cm

\bibliographystyle{amsplain}

\newpage

\vskip 1cm

2000 {\bf Mathematics Subject Classification.} 35B40, 35K55, 35K65, 35J87, 26A33.

{\bf Keywords and phases.} Porous medium equation, fractional Laplacian, obstacle problem, asymptotic behavior.


\

{\sc Addresses:}

\medskip

{\sc Luis A. Caffarelli}\newline
School of Mathematics, Univ. of Texas at Austin,
1 University Station, C1200, Austin, Texas 78712-1082. \newline
Second affiliation: Institute for Computational Engineering and Sciences.\newline
e-mail: caffarel@math.utexas.edu

\medskip

{\sc Juan Luis V{\'a}zquez}\newline
Departamento de Matem\'{a}ticas, Universidad Aut\'{o}noma de Madrid, 28049
Madrid, Spain.  Second affiliation: Institute ICMAT. \newline
e-mail: juanluis.vazquez@uam.es

\end{document}